# Shape optimization of flywheel used in agricultural thresher

Prafull Bhosale
Mechanical department
Indian Institute of Technology Bombay,
Mumbai, India
203010005@iitb.ac.in

Umesh Zawar
Mechanical department
Indian Institute of Technology Bombay,
Mumbai, India
203100044@iitb.ac.in

*Abstract*— **This research paper presents the procedure for shape optimization of flywheel used in an agricultural thresher machine using a cubic B-spline curve and the Jaya algorithm. The flywheel is an essential element for storing kinetic energy in modern machines. Shape optimization of the flywheel was carried out using a cubic B-spline curve to maximize the kinetic energy of the flywheel subjected to constraints of mass of flywheel and permissible value of Von-Mises stress in a flywheel. The B-spline curve was used as it gives local control and is not restricted by control points. The control points of the B-spline are considered design variables. The Von-Mises stresses at all points in the flywheel are determined by solving two-point boundary value differential equations using the finite difference method. Finally, the optimization problem was solved using the Jaya algorithm. These results are then compared with results obtained using other nature-inspired algorithms like a genetic algorithm (GA) and particle swarm optimization (PSO), and it is found that the Jaya algorithm gives better results. This project is based on the research paper optimal design of the flywheel using nature-inspired optimization algorithms (Prem Singh, Himanshu Choudhary), which applied GA, PSO, and Jaya algorithms. [1]**

The developed MATLAB code can be downloaded from https://github.com/prafull-bhosale/Flywheel_Shape_Optim

## I. INTRODUCTION

The flywheel is an energy storage device and helps to smooth the torque fluctuations by storing and releasing energy. Generally, it is used in many applications such as internal combustion engines, punch presses, agricultural thresher machines, vehicle power stations, aerospace, etc., due to its high efficiency, high cyclic lifetime, easy maintenance, and low environmental pollution [16]. The performance of the flywheel can be increased by increasing its energy storing capacity. The energy storing capacity of the flywheel can be determined by its material property, rotating speed, and geometry. A material with high density can significantly increase the flywheel mass and thus increase its energy storing capacity. Further, the high-strength material can withstand more stress. As a result, the flywheel can rotate at a higher speed and store more kinetic energy. Unfortunately, the number of these materials that can be chosen for most designs is limited due to their availability. Moreover, a large mass flywheel is restricted by the available space in the machine and can also transfer large gravitational forces to the bearings. Fortunately, the geometry of a flywheel can also influence the storage capacity of kinetic energy for the same flywheel mass. Composite materials [17],[18],[19],[20] are widely used to achieve long life and high performance in the design of the geometry of a flywheel. But, the cost of composite materials is high compared to metallic materials, and the forming process is also relatively complex. Therefore, metallic materials are used to design the various geometries of flywheels due to their low price and simple manufacturing process. Recently, most studies have been focused on shape optimization of the flywheel rotor cross-section. It is an efficient method to increase the kinetic energy by optimizing the rotor thickness along a radial direction [21].

In another study, coefficients of a Fourier series and Fourier sine series are considered as design variables. The ratio of inertia to volume and kinetic energy are considered objective functions with constraints of the maximum allowable stress, maximum thickness, and maximum mass for the optimal shape of the flywheel. These series represent the thickness as a function of the radius. The radial and tangential stresses of the flywheel are computed using a two-point boundary value differential equation [7][8]. Further, a non-linear optimization problem is formulated to maximize the energy density (stored energy per unit mass) using the parametric geometry modeling method. Then, the downhill simplex method was used to solve this problem [21]

Unfortunately, the polynomial expansion, Fourier series, and Fourier sine series are limited by the number of coefficients and do not describe the degree of the resulting curve. Although Bezier curves describe the degree of the curve, these curves are controlled globally and are limited by control points. Moreover, the optimal shape of the flywheel is determined using conventional optimization techniques and finite element analysis. But these methods increase the computational cost and are also less efficient. In this paper, the shape optimization model of the flywheel is formulated using a cubic B-spline curve. The B-Spline curves have local control and are not limited by control points. Then, the optimization problem is solved by the Jaya algorithm. The effectiveness of the approach is illustrated by the flywheel design of an agricultural thresher machine.

The thresher is a piece of agricultural machinery that detaches the grains from harvested crops by a combination of impact and a rubbing action with minimum effort and time. Fig. 1 The flywheel is an essential element of the thresher machine [22]. The function of the flywheel in a thresher machine is to minimize the variations in the speed of the PTO

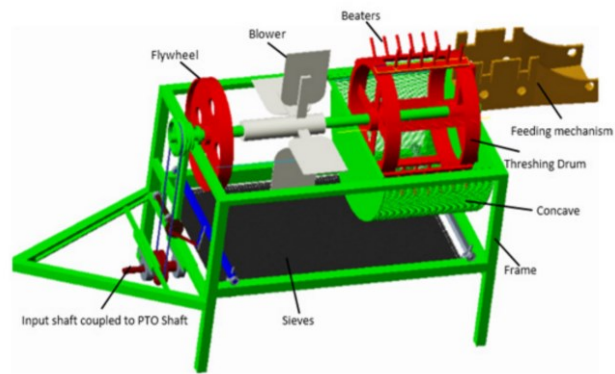

Fig. 1: Agricultural thresher machine with flywheel

shaft due to torque fluctuations of the threshing drum by storing or releasing kinetic energy.

## II. MATERIAL AND METHODS

### A. Cubic B-spline curve

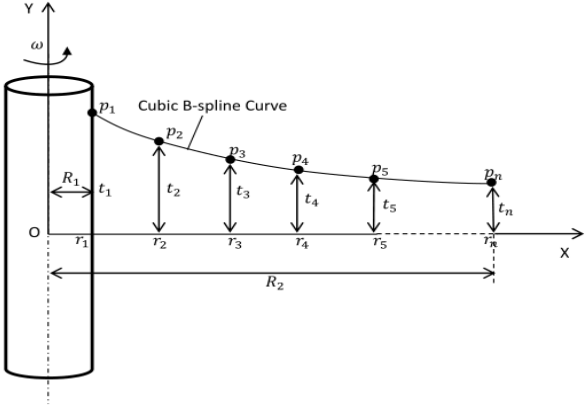

Fig. 2: Cross-section profile of flywheel

The cross-section profile of the flywheel is represented by a cubic B-spline as shown in Fig. 2. R1 and R2 are the inner and outer radii of the flywheel, respectively. The flywheel is axisymmetric about the Y axis and also symmetric about the X axis. The thickness of the flywheel along the radial direction is represented by a cubic B-spline. [10][11] A set of control points defined the actual curve as given by parametric equation (1). X and Y coordinates of control points are taken as radius and thickness of flywheel, respectively.

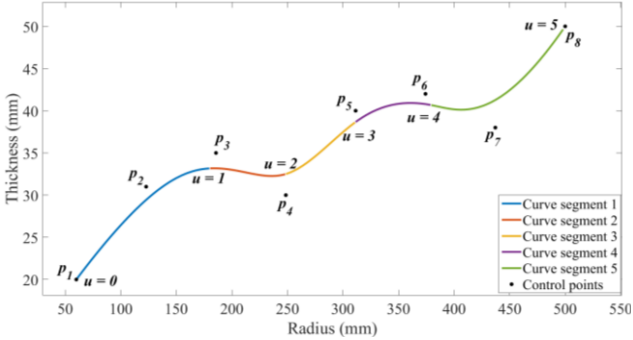

Fig. 3: Nonperiodic B-spline curve segments with control points

A nonperiodic cubic B-spline made up of curve segments connected with $C^2$ continuity, as shown in Fig. 3.

$$\boldsymbol{p}(u) = \sum_{i=1}^{n} N_{i,k}(u) * \boldsymbol{p}_i, \qquad 0 \leq u \leq S \quad (1)$$

where parameter $k$ controls the degree of the curve. $k = 4$ for cubic B-spline.
$N_{i,k}(u)$ are the B-spline basis functions, which are defined as follows:

$$N_{i,k} = (u - v_i)\frac{N_{i,k-1}(u)}{v_{i+k-1} - v_i} + (v_{i+k} - u)\frac{N_{i+1,k-1}(u)}{v_{i+k} - v_{i+1}} \quad (2)$$

where

$$N_{i,1} = \begin{cases} 1, & v_i \leq u \leq v_{i+1} \\ 0, & \text{elsewhere} \end{cases} \quad (3)$$

$v_i$ are the knot values that decide which control points contribute to defining a particular curve segment. These values depend on whether the B-spline is open (nonperiodic) or periodic curve. The periodic B-spline curve does not pass through control points, whereas the nonperiodic curve passes through two endpoints.

The knot values for the nonperiodic B-spline curve are as follows: [12]

$$v_i = \begin{cases} 0, & i < k \\ i - k + 1, & k \leq i \leq n \\ n - k + 1, & i > n \end{cases} \quad (4)$$

where $\qquad 0 \leq i \leq n + k - 1$

The coordinates of any point as a function of $u$ on the $i^{th}$ the segment of the curve for periodic B-spline curves is determined using (1), (2), (3), and (4) and are given as

$$\boldsymbol{p}^{(i)}(u) = a_1^{(i)}\boldsymbol{p}_i + a_2^{(i)}\boldsymbol{p}_{i+1} + a_3^{(i)}\boldsymbol{p}_{i+2} + a_4^{(i)}\boldsymbol{p}_{i+3} \quad (5)$$

where
$$\boldsymbol{p}^{(i)}(u) = \begin{bmatrix} r^{(i)}(u) & t^{(i)}(u) \end{bmatrix}^T$$
$$\boldsymbol{p}_i = [r_i \quad t_i]^T$$
$$\boldsymbol{p}_{i+1} = [r_{i+1} \quad t_{i+1}]^T$$
$$\boldsymbol{p}_{i+2} = [r_{i+2} \quad t_{i+2}]^T$$
$$\boldsymbol{p}_{i+3} = [r_{i+3} \quad t_{i+3}]^T$$

$\boldsymbol{p}_i, \boldsymbol{p}_{i+1}, \boldsymbol{p}_{i+2}, \boldsymbol{p}_{i+3}$ are the control points and

$a_1^{(i)}, a_2^{(i)}, a_3^{(i)}, a_4^{(i)}$ are obtained by solving (2), (3), and (5) for each segment in MATLAB.

Thus, the equations for each segment become:

$$\begin{aligned} \boldsymbol{p}^{(1)}(u) = &-\boldsymbol{p}_1 (u - 1)^3 \\ &+ 0.25\, \boldsymbol{p}_2\, u\, (7\, u^2 - 18\, u + 12) \\ &- 0.0833\, \boldsymbol{p}_3\, u^2\, (11\, u - 18) \\ &+ 0.1667\, \boldsymbol{p}_4\, u^3 \end{aligned} \quad (6)$$

$$\begin{aligned} \boldsymbol{p}^{(2)}(u) = &-0.25\, \boldsymbol{p}_2\, (u - 2)^3 \\ &+ \boldsymbol{p}_3\, (0.5833\, u^3 - 3\, u^2 + 4.5\, u - 1.5) \\ &- \boldsymbol{p}_4\, (0.5\, u^3 - 2\, u^2 + 2\, u - 0.6667) \\ &+ 0.1667\, \boldsymbol{p}_5\, (u - 1)^3 \end{aligned} \quad (7)$$

$$\begin{aligned} \boldsymbol{p}^{(3)}(u) = &-0.1667\, \boldsymbol{p}_3\, (u - 3)^3 \\ &+ \boldsymbol{p}_4\, (0.5\, u^3 - 4\, u^2 + 10\, u - 7.3333) \\ &- \boldsymbol{p}_5\, (0.5\, u^3 - 3.5\, u^2 + 7.5\, u - 5.1667) + 0.1667\, \boldsymbol{p}_6\, (u - 2)^3 \end{aligned} \quad (8)$$

$$\begin{aligned} \boldsymbol{p}^{(4)}(u) = &-0.1667\, \boldsymbol{p}_4\, (u - 4)^3 \\ &+ \boldsymbol{p}_5\, (0.5\, u^3 - 5.5\, u^2 + 19.5\, u - 21.8333) \\ &- \boldsymbol{p}_6\, (0.5833\, u^3 - 5.75\, u^2 + 18.25\, u - 18.9167) + 0.25\, \boldsymbol{p}_7\, (u - 3)^3 \end{aligned} \quad (9)$$

$$\begin{aligned} \boldsymbol{p}^{(5)}(u) = &-0.1667\, \boldsymbol{p}_5\, (u - 5)^3 \\ &+ 0.0833\, \boldsymbol{p}_6\, (11\, u - 37)\, (u - 5)^2 \\ &- \boldsymbol{p}_7\, (1.75\, u^3 - 21.75\, u^2 + 89.25\, u - 121.25) + \boldsymbol{p}_8 (u - 4)^3 \end{aligned} \quad (10)$$

## B. Flywheel parameters

The mass and kinetic energy are parameters of the flywheel which describe its performance. These parameters for each segment are calculated as follows:

$$m^{(i)} = \rho \int_{u_{i-1}}^{u_i} dV \quad (11)$$

$$m^{(i)} = 2\pi\rho \int_{u_{i-1}}^{u_i} t^{(i)}(u) r^{(i)}(u) \frac{dr^{(i)}(u)}{du} du \quad (12)$$

$$e_k^{(i)} = \frac{1}{2}\rho\omega^2 \int_{u_{i-1}}^{u_i} \left(r^{(i)}(u)\right)^2 dV \quad (13)$$

$$e_k^{(i)} = \pi\rho\omega^2 \int_{u_{i-1}}^{u_i} t^{(i)}(u) \left(r^{(i)}(u)\right)^3 \frac{dr^{(i)}(u)}{du} du \quad (14)$$

where $\rho$ is the density of the flywheel material and $\omega$ represents the angular velocity of the flywheel. Further, the whole mass and kinetic energy of the flywheel are calculated by the summing of the mass, and kinetic energy of each segment, respectively as:

$$\begin{aligned} M = \sum_{i=1}^{S} m^{(i)} &= 157.5454\, t_1 + 334.8153\, t_2 \\ &+ 500.4834\, t_3 + 733.6096\, t_4 \\ &+ 942.2381\, t_5 + 1.0172e+03\, t_6 \\ &+ 1.0883e+03\, t_7 + 837.9416\, t_8 \end{aligned} \quad (15)$$

$$\begin{aligned} E_k = \sum_{i=1}^{S} e_k^{(i)} &= 3.1720e+03\, t_1 + 1.4887e+04\, t_2 \\ &+ 4.1342e+04\, t_3 + 1.0174e \\ &+ 05\, t_4 + 2.0967e+05\, t_5 \\ &+ 3.1681e+05\, t_6 + 4.3521e \\ &+ 05\, t_7 + 4.0134e+05\, t_8 \end{aligned} \quad (16)$$

## C. Stress analysis of flywheel

A nonperiodic uniform B-spline cubic curve is used for stress analysis of the flywheel. It passes through the endpoints. Stresses at the inner and outer radii are calculated using this curve. While a periodic uniform B-spline cubic curve does not pass through the endpoints, it does not give information about stresses at the inner and outer radius of the flywheel.

Tangential and radial stresses occur due to centrifugal forces during flywheel operation. The stress distribution in a flywheel of arbitrary shape is accomplished by applying force balance on a small element of the flywheel. [4], [6] Such an element is shown in Fig. 4. Assuming stresses due to tangential forces are small compared to centrifugal forces, the force equilibrium equation for the element is as follows:

$$\frac{d}{dr^i(u)}(t^i(u) r^i(u) \sigma_r^i) - t^i(u)\sigma_\theta^i + \rho(r^i(u))^2 t^i(u)\omega^2 = 0 \quad (17)$$

$$(\sigma_\theta^{(i)} - \sigma_r^{(i)})(1+v) + r\frac{d\sigma_\theta^{(i)}}{dr^{(i)}(u)} - r^{(i)}(u)v\frac{d\sigma_r^{(i)}}{dr^{(i)}(u)} = 0 \quad (18)$$

where $\sigma_r^i$ and $\sigma_\theta^i$ are tangential and radial stresses for each segment, respectively. $t^i(u)$ and $r^i(u)$ are thickness and radius of the flywheel for each segment given in (5)

Let us define a stress function for each segment as

$$Z^{(i)} = t^{(i)}(u) r^{(i)}(u) \sigma_r^{(i)} \quad (19)$$

$$\sigma_r^{(i)} = \frac{Z^{(i)}}{t^{(i)}(u) r^{(i)}(u)} \quad (20)$$

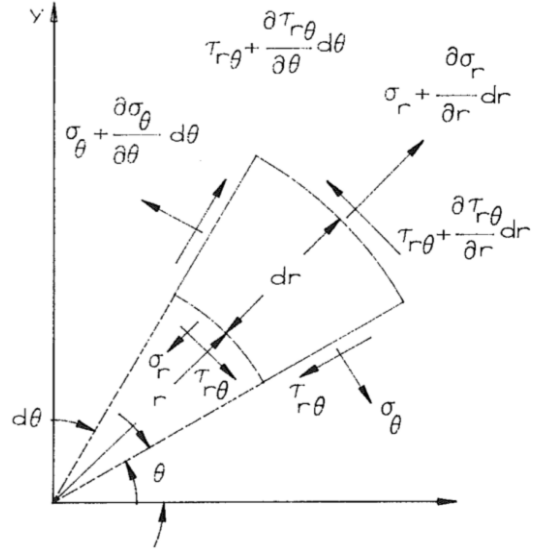

Fig. 4: General element for stress derivation in cylindrical coordinates

Substituting (20) in (17) and solving for $\sigma_\theta$ we get:

$$\sigma_\theta^{(i)} = \frac{1}{t^{(i)}(u)} \left( \frac{dZ^{(i)}}{dr^{(i)}(u)} + \rho(r^{(i)}(u))^2 \omega^2\, t^{(i)}(u) \right) \quad (21)$$

A second-order ordinary differential equation is obtained by solving (17), (18), and (20) as:

$$r^2 \frac{d^2 Z}{dr^2} + r \frac{dZ}{dr} - Z + (3+v)\rho r^3 \omega^2 t - \frac{r}{t}\frac{dt}{dr}\left(r\frac{dZ}{dr} - vZ\right) = 0 \quad (22)$$

Equations (21) and (22) are written in parametric form for each segment by converting the independent variable $r$ into $u$ by using chain rule of differentiation as:

$$\sigma_\theta^{(i)} = \frac{1}{t^{(i)}(u)} \left( \frac{\frac{dZ^{(i)}}{du}}{\frac{dr^{(i)}(u)}{du}} + \rho(r^{(i)}(u))^2 \omega^2\, t^{(i)}(u) \right) \quad (23)$$

$$C^{(i)}(u) \frac{d^2 Z^{(i)}}{du^2} + D^{(i)}(u) \frac{dZ^{(i)}}{du} + E^{(i)}(u) Z^{(i)} = F^{(i)}(u) \quad (24)$$

where

$$C^{(i)}(u) = \left(r^{(i)}(u)\right)^2 \frac{dr^{(i)}(u)}{du} \quad (25)$$

$$D^{(i)}(u) = r^{(i)}(u)\left(\frac{dr^{(i)}(u)}{du}\right)^2 - \left(r^{(i)}(u)\right)^2 \frac{d^2 r^{(i)}(u)}{du^2}$$
$$- \frac{\left(r^{(i)}(u)\right)^2}{t^{(i)}(u)} \frac{dr^{(i)}(u)}{du} \frac{dt^{(i)}(u)}{du} \quad (26)$$

$$E^{(i)}(u) = v \frac{r^{(i)}(u)}{t^{(i)}(u)} \frac{dt^{(i)}(u)}{du}\left(\frac{dr^{(i)}(u)}{du}\right)^2 - \left(\frac{dr^{(i)}(u)}{du}\right)^3 \quad (27)$$

$$F^{(i)}(u) = -(3+v)\rho\omega^2 t^{(i)}(u)\left(r^{(i)}(u)\right)^3 \left(\frac{dr^{(i)}(u)}{du}\right)^3 \quad (28)$$

Equation (22) is a non-linear second-order ordinary differential equation with independent variable $u$ and dependent variable $Z$. For flywheel, we have two boundary conditions; radial stress at inner and outer radius is zero. So, stress function Z is also zero at the inner and outer radius (20).

Then (22) becomes a two-point boundary value problem with second order non-linear differential equation. It can be solved by ode45 and using the Runge-Kutta method combined with the shooting method. [1],[15].

In the finite difference method, we approximate the first and second derivative of functions as follows: [24]
Central difference approximation for the first and second derivative of Z with respect to $u$ is given by (29) and (30)

$$DZ_j = \frac{Z_{j+1} - Z_{j-1}}{2h} \quad (29)$$

$$D^2 Z_j = \frac{Z_{j+1} - 2Z_j + Z_{j-1}}{h^2} \quad (30)$$

where $D_u \equiv \frac{d}{du}$, $D_u^2 \equiv \frac{d^2}{du^2}$, $h = u_j - u_{j-1}$

These approximations have errors in the order of $h^2$. These are generally the most accurate approximations, but in some situations, for example, at the left and right boundary nodes, it is necessary to use forward or backward difference approximations.

e.g., the forward difference approximation for $DZ_j$ is:
$$DZ_j = \frac{Z_{j+1} - Z_j}{h} \quad (31)$$

The backward difference approximation for $Dz_j$ is:
$$DZ_j = \frac{Z_j - Z_{j-1}}{h} \quad (32)$$

After this, we get $n$ equations for $n$ nodes in $n$ unknown $Z_j$ which are then solved simultaneously to get values of stress function at each node.

Once the stress function is known at various discrete points in each segment, we can calculate radial and tangential stresses using (20) and (21). Then we can determine Von-Mises at each point using distortion energy theory as:

$$\sigma_t = (\sigma_r^2 + \sigma_\theta^2 - \sigma_r \sigma_\theta)^{1/2} \quad (33)$$

This allows us to satisfy stress constraints for the safe design of the flywheel.

### III. FORMULATION OF AN OPTIMIZATION PROBLEM

A non-linear constrained optimization problem for optimal thickness distribution along the radial direction is formulated using the cubic B-spline. The design parameters of the flywheel are shown in Fig. 4. OXY represents the global coordinate system. The coordinates of control points $r_1\ r_2\ r_3\ r_4 \ldots\ldots\ldots r_n$ in X and $t_1\ t_2\ t_3\ t_4 \ldots\ldots\ldots t_n$ in Y direction are radius and the thickness of the flywheel, respectively.

#### A. Design variables

The radius of the flywheel between $R_1$ and $R_2$ is equally divided. Thus, the X-coordinate of control points is fixed given as follows:

$$r_i = R_1 + \frac{(R_2 - R_1)(i-1)}{n-1} \quad (34)$$

where $i=1, 2, 3, \ldots, n$

The Y-coordinates of control points are thicknesses that are considered as design variables and are expressed in vector form as:

$$\mathbf{x} = [t_1\ t_2\ t_3\ t_4 \ldots\ldots\ldots t_n]^T \quad (35)$$

Where $n$ = Number of control points

#### B. Objective function and constraints

For better performance, the flywheel should store as much kinetic energy as possible. The kinetic energy of the flywheel is determined by its shape. Thus, the optimization problem can be stated as 'To find the optimal shape of the flywheel, which maximizes the kinetic energy.' The kinetic energy of the flywheel is determined by using (16).

Here we want to maximize the kinetic energy. So, the minimization problem can be stated as follows:
$$\text{Minimize } f(\mathbf{x}) = -E_k \quad (36)$$

The mass of the flywheel should be less than the maximum allowable mass for a given application ($M_{max}$). The mass of the flywheel can be determined using (15). The mass constraint is defined as:

$$g_1(\mathbf{x}) = M \leq M_{max} \quad (37)$$

Further, the maximum value of Von-Mises stress at any point in the flywheel must be less than permissible stress $\sigma_a$. The maximum value of Von Mises stress at all points of each segment in the radial direction is calculated using (33). The stress constraint is defined as:

$$g_2(\mathbf{x}) = max(\sigma_t) \leq \sigma_a \quad (38)$$

Finally, the optimization problem is posed as:

$$\text{Minimize } f(\mathbf{x}) = -E_k$$

subjected to
$$g_1(\mathbf{x}) = M - M_{max} \leq 0 \quad (39)$$
$$g_2(\mathbf{x}) = max(\sigma_t) - \sigma_a \leq 0$$
$$LB_i \leq x_i \leq UB_i \quad i = 1,\ldots,n$$

The negative sign in the expression for the objective function implies that the function $E_k$ is being maximized. $LB_i$ and $UB_i$ are the lower and upper bounds on the $i^{th}$ design variable, and $n$ represents the number of design variables. To obtain an optimum solution, the constrained problem, as defined in (39), is converted into an unconstrained problem by using a penalty function [5]. A significant penalty value is added to the objective function for each constraint violation. As a result, the objective function proceeds toward an infeasible solution. Hence, the global optimum solution is obtained by satisfying all the constraints using a suitable optimization algorithm. The original constrained optimization problem is then stated as an unconstrained optimization problem in which the first and second term describes the objective function and the penalty function, respectively. Finally, the shape optimization problem of the flywheel is formulated as:

$$Minimize\ f(x) = -E_k + \sum_{B=1}^{2} C_B(CP)^B \quad (40)$$

The upper and lower limit of variables are:

and $r1_{j,i}$ and $r2_{j,i}$ are the two random numbers for the $j^{th}$ variable during the $i^{th}$ iteration in the range [0, 1].

The term "$r1_{j,i}(x_{j,best,i} - |x_{j,k,i}|)$" indicates the tendency of the solution to move closer to the best solution, and the term "$-r2_{j,i}(x_{j,worst,i} - |x_{j,k,i}|)$" indicates the tendency of the solution to avoid the worst solution. $x'_{j,k,i}$ is accepted if it gives a better function value. All the accepted function values at the end of the iteration are maintained, and these values become the input to the next iteration. The algorithm always

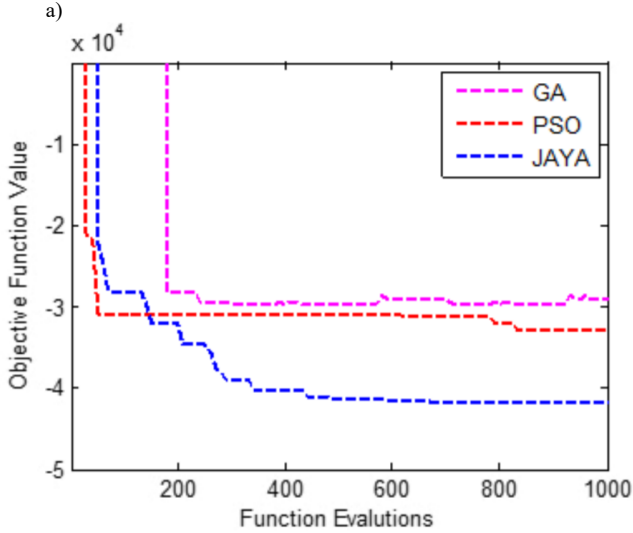
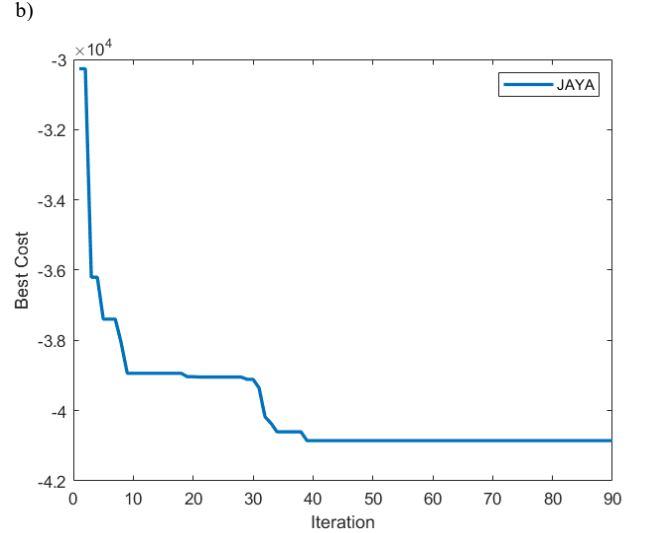

Fig. 5: Convergence plots of best objective function value in (a) research paper vs (b) Jaya algorithm for our project

$$LB_i \leq x_i \leq UB_i \quad i = 1,..,n \quad (41)$$

where $CP$ is the penalty value of the order of $10^8$, which is added to the objective function if the constraints are not satisfied and $C_B$ is the Boolean function defined as [5]

$$C_B = \begin{cases} 0 & if\ g_B(x) \leq 0 \\ 1 & otherwise \end{cases} \quad (42)$$

C. Jaya algorithm

Jaya algorithm is based on the concept that the solution obtained for a given problem should move towards the best solution and should avoid the worst solution [14]. This algorithm requires only the common control parameters and does not require any algorithm-specific control parameters. Let $f(x)$ is the objective function to be minimized (or maximized). At any iteration $i$, assume that there are '$m$' number of design variables (i.e., $j=1, 2, ..., m$), '$n$' number of candidate solutions (i.e., population size, $k=1, 2..., n$). Let the best candidate obtains the best value of $f(x)$ (i.e., $f(x)_{best}$) in the entire candidate solutions, and the worst candidate obtain the worst value of $f(x)$ (i.e., $f(x)_{worst}$) in the entire candidate solutions. If $x_{j,k,i}$ is the value of the $j^{th}$ variable for the $k^{th}$ candidate during the $i^{th}$ iteration, then this value is modified as per the following (43).

$$x'_{j,k,i} = x_{j,k,i} + r1_{j,i}(x_{j,best,i} - |x_{j,k,i}|) \\ - r2_{j,i}(x_{j,worst,i} - |x_{j,k,i}|) \quad (43)$$

where $x_{j,best,i}$ is the value of the variable $j$ for the best candidate and $x_{j,worst,i}$ is the value of the variable $j$ for the worst candidate. $x'_{j,k,i}$ is the updated value of $x_{j,k,i}$ tries to get closer to success (i.e., reaching the best solution) and tries to avoid failure (i.e., moving away from the worst solution). The algorithm strives to become victorious by getting the best solution, and hence it is named as Jaya (a Sanskrit word meaning victory).

Notations used in the Jaya algorithm are defined as follows:
pop = population size = 1000
$LB_i$, $UB_i$ = Lower and upper bounds for the $i^{th}$ design variables
$x_i$ = $i^{th}$ design variable
$x_{ij}$ = $i^{th}$ design variable for $j^{th}$ population
$f_j$ = An objective function value of the $j_{th}$ population
rand = Any random number in the range of 0 and 1
MaxStallGenerations = 50
FunctionTolerance = $10^{-6}$

After formulating the optimization problem, it can be solved either by classical methods or nature-inspired optimization algorithms. The traditional methods involve the calculation of the gradient of the objective function, which increases computational complexity if the problem is non-linear. Nature Inspired Optimization Algorithms have less computational complexity compared to classical methods. There are some nature-inspired optimization algorithms like Genetic Algorithm (GA), Particle Swarm Optimization (PSO), and Ant Colony Optimization (ACO), which converge to a global minimum, but there is no guarantee of an optimal global solution. Some of these algorithms, such as GA, PSO, and ACO, have algorithmic control parameters, which affect the performance of algorithms. The Jaya algorithm is a parameter-less algorithm that has only one phase, unlike the two phases of the TLBO algorithm [13]. It gives the optimal solution quickly and updates the worst solution [14]. A flow chart of the Jaya algorithm for the shape of the flywheel is shown in Fig. 8.

Here, the stopping criterion for the Jaya algorithm is considered as the number of iterations and/or maximum stall generations. Maximum stall generations are the number of generations for which the difference between consecutive optimization function values is less than some tolerance. [3]

$$\frac{|f_j - f_{j-1}|}{max(1, f_j)} \leq \text{FunctionTolerance} \qquad (44)$$

where $f_j$ is the value of the objective function for $j^{th}$ generation and $f_{j-1}$ is the value of the objective function for $j - 1^{th}$ generation.

The algorithm will stop if the maximum stall generations exceed the user-specified value. Therefore, design variables do not affect the function evaluations but may affect the computational time of the algorithm.

## IV. RESULTS AND DISCUSSION

The above optimization procedure is applied to the flywheel of an agricultural thresher machine. The material properties and design parameters of the flywheel as given in TABLE I. and TABLE II. respectively. The design variables are the Y-coordinates of the control points as discussed in III.A. The lower and upper bounds of the design variables are taken as:

$$0.010 \leq x_i \leq 0.06 \qquad i = 1, 2, \ldots, 8 \qquad (45)$$

The second-order differential equation for stress function is solved using the finite difference method using a step size of $h=0.01$. Further, the radial, tangential, and Von-Mises stresses are calculated using (20), (21), and (33), respectively. The stress distribution for optimum values of design variables is shown in Fig. 6.

The Jaya algorithm was coded in MATLAB. The algorithm was run multiple times to make sure that the function attained the global minimum. The convergence plot for one of such runs is shown in Fig. 5. For comparing the efficiency of the Jaya algorithm, the results obtained by using the Jaya algorithm are compared with the GA and PSO results as shown in 0. It was found that the Jaya algorithm gives better results compared to GA and PSO.

TABLE I. MATERIAL PROPERTIES OF FLYWHEEL

| Material | Density(kg/$m^3$) | Elastic Modulus (GPa) | Poisson's ratio |
|---|---|---|---|
| Grey cast iron | 7250 | 210 | 0.3 |

TABLE II. DESIGN PARAMETERS OF THE FLYWHEEL

| Control points | Inner Radius $R_1$ (mm) | Outer Radius $R_2$(mm) | Angular Velocity $\omega$ (rad/s) | $M_{max}$ (kg) | $\sigma_a$ ($N/mm^2$) |
|---|---|---|---|---|---|
| 8 | 0.06 | 0.5 | 65.45 | 115 | 6.4 |

The optimized shape of the flywheel is shown in Fig. 7. We can see that the Jaya algorithm tries to move as much material as possible to the outer radius of the flywheel so as to increase its kinetic energy storage for the given mass and stress constraints. However, more the material situated at the outer part of the flywheel, the centrifugal force acting on material close to the inner radius. So, it is evident from the figure that the algorithm also tries to keep some material at the inner radius so that Von-Mises stresses are less than allowable stress and the stress constraint is satisfied.

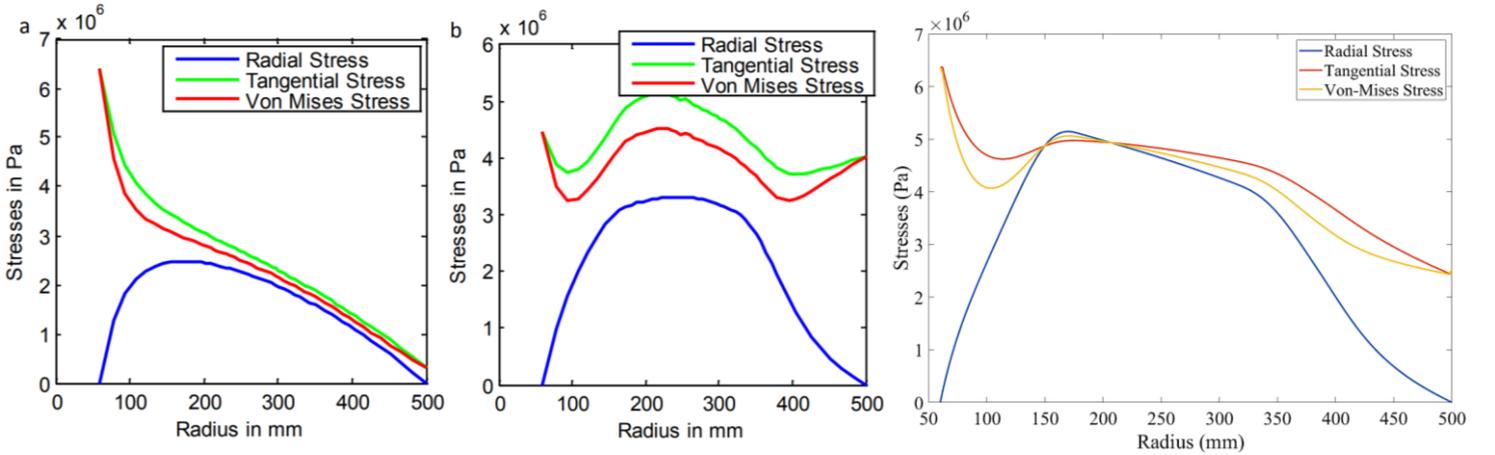
Fig. 6: Stress distribution along radial direction in a) original flywheel of constant thickness b) optimized flywheel research paper c) optimized flywheel our project

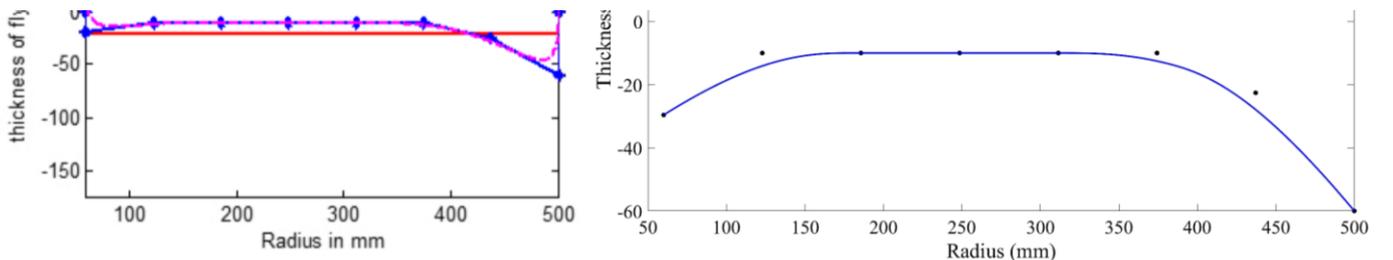
Fig. 7: Optimized shape of flywheel a) research paper b) our project

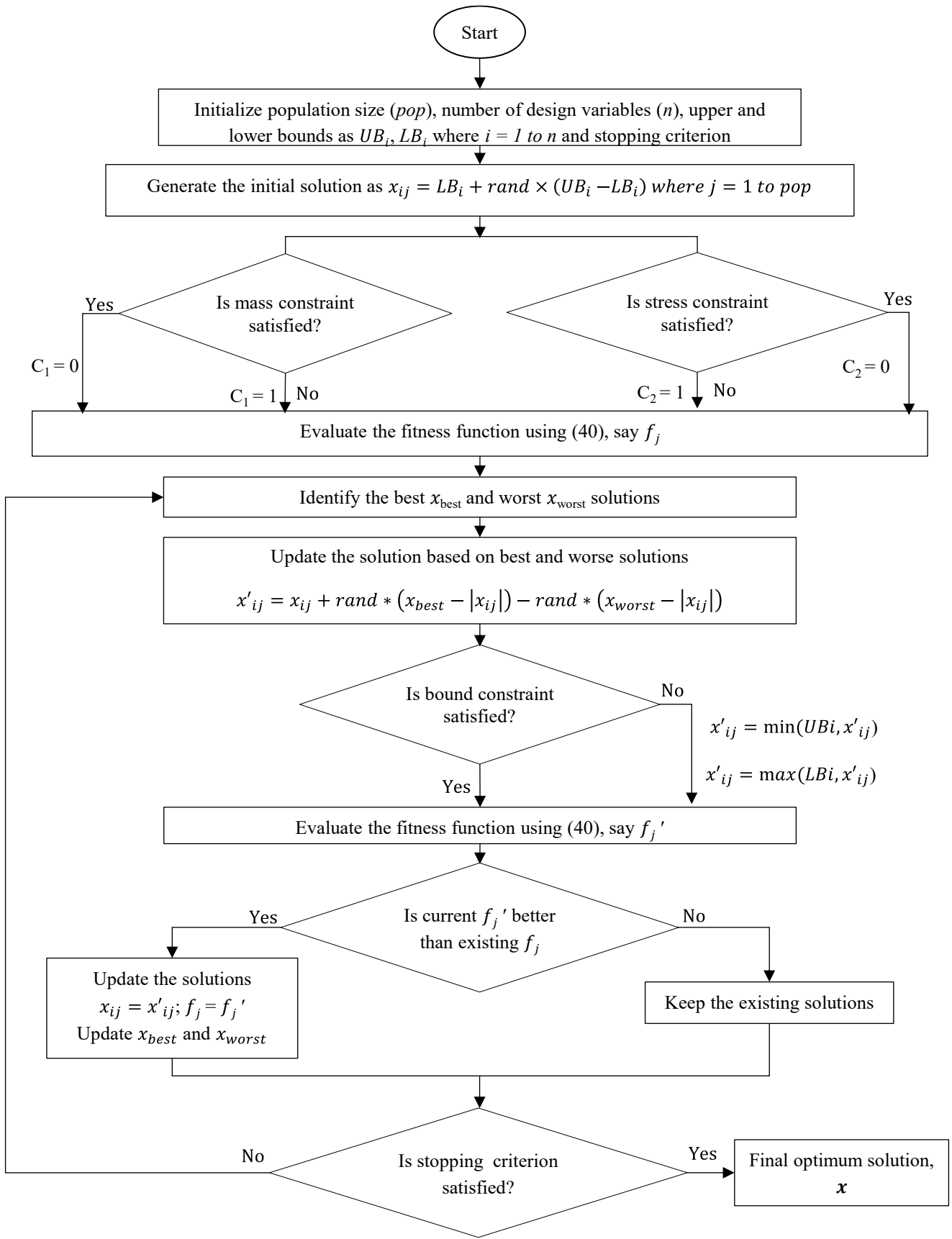

Fig. 8: Jaya algorithm flowchart for shape optimization of flywheel

TABLE III. OPTIMUM VALUES FOR THE SHAPE OF THE FLYWHEEL

| Optimum Values | Constant thickness flywheel | GA | PSO | Jaya from research paper | Jaya from our project |
|---|---|---|---|---|---|
| $t_1(m)$ | 0.020 | 0.0267 | 0.01 | 0.0189 | 0.0296 |
| $t_2(m)$ | 0.020 | 0.0122 | 0.01 | 0.01 | 0.01 |
| $t_3(m)$ | 0.020 | 0.0228 | 0.06 | 0.01 | 0.01 |
| $t_4(m)$ | 0.020 | 0.0292 | 0.01 | 0.01 | 0.01 |
| $t_5(m)$ | 0.020 | 0.0248 | 0.01 | 0.01 | 0.01 |
| $t_6(m)$ | 0.020 | 0.0182 | 0.01 | 0.01 | 0.01 |
| $t_7(m)$ | 0.020 | 0.0013 | 0.01 | 0.0246 | 0.0226 |
| $t_8(m)$ | 0.020 | 0.0211 | 0.0489 | 0.06 | 0.06 |
| Mass *(kg)* | 115 | 115 | 115 | 115 | 114.96 |
| Kinetic Energy *(J)* | 30483.66 | 28950 | 32896.25 | **41624.88** | **40854.16** |
| *f(x)* | - | -28950 | -32896.25 | **-41624.88** | **-40854.16** |

## V. CONCLUSION

The procedure for the optimal design of flywheel using cubic B-spline and Jaya algorithm was described in this paper. The optimization problem was formulated to maximize the kinetic energy stored in the flywheel subjected to mass and Von-Mises stress constraints. The mentioned approach is applied to the flywheel of an agricultural thresher machine. It was found that the algorithm tries to concentrate as much material as possible to the outer radius of the flywheel so as to maximize the kinetic energy storage. This is as per our intuition. This creates higher stress levels at interior locations in the flywheel. This trade-off between moving material toward the outside radius and keeping the stress within acceptable limits is evident in Fig. 6 and Fig. 7. Thus, the profile of the optimized flywheel can be divided into three sections; a thin middle section and a thicker inner and outer section. It was found that the Jaya algorithm gives better results compared to GA and PSO. The optimized shape of the flywheel, obtained by using by Jaya algorithm, stores 34.02 % more energy compared with that of the original flywheel, without failure during its operation. As a result of this, the dynamic performance (like vibrations and torque fluctuation) of the existing thresher machine will improve.

It is found that the Jaya algorithm gives better results compared to GA and PSO.